\documentclass{article}
\newcommand{\bs}[1]{\boldsymbol{#1}}
\usepackage{authblk}
\usepackage{amsmath}
\usepackage{amssymb}
\usepackage[pdftex]{graphicx}
\usepackage{color}

\usepackage[margin=20truemm]{geometry}
\title{\textcolor{black}{Automatically well-conditioned} \textcolor{black}{collocation boundary element method for transmission problems based on the Burton--Miller formulation}}
\author[1]{Yasuhiro Matsumoto$^*$}
\author[2]{Akihiro Yoshiki$^*$}
\author[2]{Hiroshi Isakari$^\dagger$}
\affil[1]{\small \textcolor{black}{Center for Information Infrastructure, Institute of Science Tokyo, Tokyo, Japan.}}
\affil[2]{\small Faculty of Science and Technology, Keio University, Kanagawa, Japan.}
\affil[*]{Contributed equally to this work.}
\affil[$\dagger$]{Corresponding author: isakari@sd.keio.ac.jp}
\date{ }

\begin{document}

\maketitle

\abstract{
This paper proposes a collocation boundary element method based on the Burton--Miller method for solving transmission problems, which is rapidly convergent within the Krylov subspace solver framework. Our study enhances Burton--Miller-type boundary integral equations tailored for transmission problems by exploiting the Calderon formula. In cases where a single material exists in an unbounded host medium, we demonstrate the formulation of the boundary integral equation such that the underlying integral operator ${\cal A}$ is spectrally well-conditioned. Specifically, ${\cal A}$ can be designed such that ${\cal A}^2$ has only a single eigenvalue accumulation point. Furthermore, we extend this to the multi-material case, proving that the square of the proposed operator has only a few eigenvalues clustering points. When the collocation method is used to discretise the proposed boundary integral equations, the good spectral properties of the integral operator are naturally inherited to the coefficient matrix $\mathsf{A}$ similarly to the Nyst\"om methods; Almost all eigenvalues of $\mathsf{A}^2$ cluster at a few points in the complex plane ensuring the small condition number for $\mathsf{A}$. Through numerical examples of several benchmark problems, we illustrate that our formulation reduces the iteration number required by iterative linear solvers, even in the presence of material junction points.}

\vspace*{10pt}


\section{Introduction}\label{sec:introduction}
Transmission problems for partial differential equations of the Helmholtz type are pervasive in various engineering fields. Efficient solvers for these problems are thus essential to develop a variety of wave devices. Classical applications encompass numerical analysis for antenna design~\cite{akduman2003direct} and forward and inverse scattering analyses in ultrasound and electromagnetic biomedical imaging~\cite{hammer1998single}. We can also mention some emerging applications such as photonic and phononic crystals, as well as metamaterials, exhibiting anomalous wave responses. These materials are expected to realise revolutionary wave devices such as superlens~\cite{pendry2000negative}, cloaking device~\cite{pendry2006controlling}, and more.

Given that the transmission problem is defined in an unbounded domain, it may be reasonable to use the boundary integral formulation to solve the problem rather than the solver of the domain type such as the finite element method (FEM). Recall that the boundary element method (BEM) bypasses the need for special considerations on the boundary condition at infinity \textcolor{black}{since it automatically gives solutions satisfying} the radiation condition. 

While various boundary integral formulations have been proposed, and some of them have widely been accepted~\cite{nedelec2001acoustic}, efficiently and accurately solving the transmission problem with BEM remains challenging in some cases. In deriving a boundary integral formulation, we need first to care about the so-called fictitious eigenvalue problem~\cite{burton1997application}. The unique solvability of a naive boundary integral equation may not be guaranteed for a specific excitation frequency even when that of the corresponding boundary value problem is. Such a frequency is called the fictitious eigenfrequency. Recent reports indicate that not only real-valued but also complex-valued fictitious eigenfrequencies may affect the accuracy of BEM~\cite{misawa2017boundary, zheng2015burton}. It is, therefore, essential to establish a boundary integral equation with fictitious eigenfrequencies as far away from the real axis as possible, considering their potential influence on accuracy. It would also be important to derive the integral equation with good spectral characteristics. We usually discretise the boundary integral equation by e.g. the Galerkin and the collocation methods and solve the resulting linear algebraic equations numerically. Since the coefficient matrix thus obtained is fully populated due to the nature of the fundamental solution, iterative solvers based on the Krylov subspace methods, such as the generalised minimal residual method~(GMRES)~\cite{saad1986gmres}, are often employed. These solvers rely on fast matrix-vector multiplication techniques such as the fast multipole method (FMM)~\cite{rokhlin1985rapid} or the hierarchical matrix method~\cite{j.ostrowski2006fast}. The number of iterations in the Krylov method depends significantly on the condition of the coefficient matrix, inherited from the spectral characteristics of the integral formulation. We thus need to carefully design the boundary integral equation.

We assert that the most widely used integral formulations for transmission problems are the PMCHWT~(Poggio--Miller--Chang--Harrington--Wu--Tsai)~\cite{chew2001fast} and M\"uller~\cite{muller1969foundations} formulations, both of which are free from real-valued fictitious eigenfrequency. The former is categorised as the integral equation of the first kind, which is thus often combined with some preconditioning techniques. Calderon's preconditioning~\cite{christiansen2000preconditionneurs, antoine2008integral} is particularly effective in improving the spectral properties of the PMCHWT equations. Niino and Nishimura pointed out that the preconditioning effect can be harnessed by a mere reordering of integral operators~\cite{niino2012preconditioning}. They also showed \textcolor{black}{that,} \textcolor{black}{upon collocation discretisation,} just reordering is more effective in reducing total computational time than multiplying the preconditioner when using GMRES. The M\"uller formulation yields a second-kind integral equation~\cite{muller1969foundations}, allowing for rapid convergence \textcolor{black}{for} the Krylov method. It is, however, noteworthy that these two integral equations may have complex-valued fictitious eigenvalue with a tiny imaginary part in its magnitude. Another possibility is the so-called single boundary integral equation (SBIE) method~\cite{costabel1985direct,kleinman1988single}. For this type of formulation of the second kind, one may augment a parameter roughly controlling the position of the complex-valued fictitious eigenvalue~\cite{misawa2016study}. The SBIE is thus considered promising in both efficiency and accuracy aspects, but applying this to multi-material cases may need some more elaboration~\cite{claeys2015secondkind}. The so-called multi-trace boundary integral formulation~\cite{hiptmair2012multiple} would also be an attractive choice, offering a straightforward extension to the multi-material problem. The formulation, however, requires multiple unknown densities on a single boundary \textcolor{black}{element. Hence it results} in slow matrix-vector multiplication in a Krylov iteration. A recently proposed formulation~\cite{vantwout2022boundary} involving both the direct and indirect boundary integrals is also efficient especially when the material constant for inclusion is considerably different from that for exterior host medium. The method may suffer from the fictitious eigenvalue problem even for real-valued frequency, though. The Burton--Miller formulation~\cite{burton1997application}, originally designed for exterior boundary value problems, can also be applied to solve the transmission problem. Although the numerical solution by this formulation with a parameter for controlling the fictitious eigenvalue is accurate, a naive use may end up with highly ill-conditioned algebraic equations. With these discussions, we may conclude that there exists no simple, efficient, accurate, and versatile boundary integral formulation for the transmission problem. 

To address the above, we here propose simple modifications to the Burton--Miller method tailored for the transmission problem. The proposed formulation

\begin{enumerate}
\item[]
\item involves a parameter to make the fictitious eigenfrequency away from the real axis in the complex plane, 
\item can easily be applied to multi-material cases even with material junctions, 
\item and has a good spectral characteristic; the square of the underlying integral operator \textcolor{black}{(implicitly computed in GMRES)} has only a few eigenvalue accumulation points.
\item[]
\end{enumerate}

\noindent The first and second features are naturally inherited by the Burton--Miller method, while the last is achieved by exploiting the Calderon formula. By referring to the formula, our formulation appropriately rearranges the order of the integral operators and slightly modifies some parts of them if necessary. By its construction, the present formulation is also benefited from 

\begin{enumerate}
\item[]
\item[4.] easy implementation. 
\item[]
\end{enumerate}

\textcolor{black}{We also mention that we have two options in discretising the boundary integral equation, i.e. the Galerkin and collocation methods. While the former fits well with theoretical discussions based on functional analysis, it requires computationally heavy numerical integrations in the practical implementation. Furthermore, the finite-dimensional matrix obtained as the Galerkin discretised integral operator does not necessarily inherit the spectral properties as is. We thus need some careful discussions with the dual basis and the Gram matrices when deriving the Calderon-preconditioned BEM; see discussions in e.g. \cite{niino2012preconditioning,hiptmair2006operator,fierro2020fast}. On the other hand, although it is somewhat inferior to the Galerkin method in terms of accuracy, the collocation BEM \textcolor{black}{does not require any special considerations during preconditioning \cite{niino2012preconditioning}. This is analogous to the fact that Calderon's preconditioning can be applied seamlessly, even when using the Nystr\"om method \cite{boubendir2016high}. Furthermore,} the collocation BEM is attractive from the perspective of simple implementation and easy use~\cite{martinsson}.  Many engineers thus employ the collocation in practical computations~\cite{isakari2012calderona,amlani2019efficient}. This is also encouraged by the fact that the hypersingular integrals in BEM can appropriately be evaluated via well-established regularisation techniques~\cite{canino1998numerical}. This paper thus discusses the collocation BEM.}

The rest of the paper is organised as follows. In Section 2, we state the boundary value problem (BVP) of interest. We then show a straightforward application of the Burton--Miller method to solve the BVP ends up with an ill-conditioned boundary integral formalism and discuss a design concept for integral formulation with good spectral characteristics in Section 3. \textcolor{black}{Section 4 reviews the well-known Calderon formula for boundary integral operators. We also numerically check that the collocation-discretised version of the operators, i.e. matrices, approximately satisfies the same relation as the Calderon formulae even when the operators are defined on the boundary with corners.} Sections \textcolor{black}{5} and \textcolor{black}{6}, the main contributions of the present paper, introduce novel Burton--Miller-type integral equations for single-material and multi-material cases, respectively. After verifying the performance of the present formulations through various numerical experiments in Section \textcolor{black}{7}, we conclude the paper and discuss future directions in Section \textcolor{black}{8}.

\section{Statement of the boundary value problem of interest}
Let us consider that $\mathbb{R}^2$ is partitioned into $M$ subdomains as $\mathbb{R}^2=\bigcup_{i=1}^{M}\Omega_i$, where $\Omega_1$ extends to infinity, while the other subdomains are all bounded. Each subdomain $\Omega_i$ equips its own constant denoted as $\varepsilon_i>0~(i=1,~\cdots,~M)$. We here normalise all the constants such that $\varepsilon_1$ is 1. This is for simplicity but does not cause any loss of generality. Given an incident field $u^\mathrm{in}:~\Omega_1\rightarrow \mathbb{C}$ satisfying the two-dimensional Helmholtz equation with the wave number $k_1>0$, the total field $u: \mathbb{R}^2\rightarrow \mathbb{C}$ solves the following BVP: 
\begin{align}
 \label{eq:helmholtz}
 \nabla^2 u(\bs{x})+k_i^2 u(\bs{x})=0~~~\bs{x}\in\Omega_i,\\
 \label{eq:jump_in_u}
 \lim_{\eta\downarrow 0}u(\bs{x}+\eta\bs{n})=\lim_{\eta\downarrow 0}u(\bs{x}-\eta\bs{n})~~~\bs{x}\in{\Gamma},\\
 \label{eq:jump_in_w}
 \lim_{\eta\downarrow 0}w(\bs{x}+\eta\bs{n})=\lim_{\eta\downarrow 0}w(\bs{x}-\eta\bs{n})~~~\bs{x}\in{\Gamma},\\
 \label{eq:radiation}
 \lim_{\|\bs{x}\|\rightarrow\infty} \|\bs{x}\|^{\frac{1}{2}} \left( \frac{\partial}{\partial\|\bs{x}\|} - \mathrm{i}k_1 \right) (u-u^\mathrm{in})(\bs{x}) = 0,
\end{align}
where $k_i:=\omega\sqrt{\varepsilon_i}$ is the wave number of the waves propagating in $\Omega_i$, and $\omega>0$ is the underlying angular frequency. Here, we assume the time dependency of all physical quantities in the time domain as $\Re[f(\bs{x})e^{-\mathrm{i}\omega t}]$, where $f$ is either $u$ or $w$, $\mathrm{i}$ is the imaginary unit, and $t$ is the time. The boundary conditions \eqref{eq:jump_in_u} and \eqref{eq:jump_in_w} on $\Gamma:=\bigcup_{i=2}^{M} \overline{\partial \Omega_i}$ define the jump relations in the total field and its flux. The flux here is defined as $w:=\frac{1}{\varepsilon_i}\nabla u \cdot \bs{n}$ with the material constants $\varepsilon_i$, where $\bs{n}$ is the unit normal on $\Gamma$. $\bs{n}$ is directed from $\Omega_1$ on $\partial\Omega_1$. The direction of $\bs{n}$ on the other boundaries is properly defined and fixed. This type of boundary value problem is often referred to as the transmission problem and is commonly found in engineering fields such as acoustics and electromagnetics. For example, a solution $u$ of the BVP \eqref{eq:helmholtz}--\eqref{eq:radiation} physically indicates the time-dependent out-of-plane component of a magnetic field, and $w$ corresponds to the in-plane tangential electric wave on the boundary. In this case, $\varepsilon_i$ represents the relative dielectric constant of a material filled in $\Omega_i$. \textcolor{black}{Note that, although we keep ourselves to two-dimensional cases in this paper, our BEM formulation in the following sections is considered effective for 3D problems almost as is.}

\section{Design concepts for the present boundary integral formalism}
Since the BVP \eqref{eq:helmholtz}--\eqref{eq:radiation} at hand is defined in an unbounded domain, the boundary integral method is a natural choice for solving it. There are various candidates for the integral formulations, and the choice may significantly affect the accuracy and efficiency of the practical computation. In this study, we ensure the accuracy by utilising the Burton--Miller method. In this section, we discuss the design concept of the integral formulation from the viewpoint of efficiency. Here, let us abstractly denote the integral equation as ${\cal A} y = b$ \textcolor{black}{corresponding to the BVP}, where $\cal A$ is an integral operator, $y$ and $b$ are respectively unknown and known density functions. The discretised version of the equation upon the collocation, i.e. the system of algebraic linear equations corresponding to the integral equation is denoted as $\mathsf{A}\bs{y}=\bs{b}$ with $\mathsf{A}$, $\bs{y}$, and $\bs{b}$ respectively being the coefficient matrix, and unknown and known vectors. \textcolor{black}{Note that the discussions presented herein are specifically geared towards the collocation method. Substantial adjustments shall be required when applying the Galerkin BEM approach~\cite{fierro2020fast}.}

\subsection{Integral equation of the second-kind}
The efficiency of the BEM largely depends on the spectral property of $\cal A$. It is thus natural to try constructing the integral equation of the Fredholm second kind, i.e. exploring $\cal A$ such that $\cal A=\cal I+\cal K$ holds, where $\cal I$ and $\cal K$ are the identity and compact operators \textcolor{black}{in an appropriate function space}, respectively. The M\"uller formulation and the SBIE are categorised into this type. It may, unfortunately, be impossible to find such a well-conditioned $\cal A$ if one employs the Burton--Miller method.

\subsection{Integral equation of the first-kind and its Calderon's preconditioning}
Another possibility may take the help of so-called preconditioning. Among the various preconditioning techniques, the Calderon preconditioning~\cite{christiansen2000preconditionneurs, antoine2008integral} is recognised as one of the best for the BEM. This method converts (or preconditions) the original integral equation as ${\cal A}{\cal M}^{-1} z = b$ such that ${\cal A}{\cal M}^{-1}={\cal I}+{\cal K}$ holds, where $\cal M$ is the right preconditioning operator, and $z:={\cal M}y$ is an auxiliary function. Niino and Nishimura~\cite{niino2012preconditioning} demonstrated that, for transmission problems, one may construct, with the PMCHWT formulation, $\cal A$ and $\cal M$ such that ${\cal M}^{-1} \simeq {\cal A}\Leftrightarrow {\cal A}^2={\cal I}+{\cal K}$ holds. \textcolor{black}{In the case that the collocation is used to discretise the integral operator ${\cal A}$, the resultant matrix $\mathsf{A}$ naturally succeeds the property; $\mathsf{A}^2$ approximately equals the identity matrix, and thus almost all eigenvalues of $\mathsf{A}^2$ exist only at a single point in the complex plane.} They also pointed out that, once the integral operator is appropriately formulated in this manner, one does not have to explicitly multiply the (inverse of) preconditioner \textcolor{black}{for the collocation BEM}. This can be understood as follows: Let us consider the preconditioned system  $\mathsf{A}^2 \bs{z}=\bs{b}$ upon discretisation and solve this by the GMRES~\cite{saad1986gmres}. Let us first recall the fact that GMRES after $n$ iterations finds the best approximated solution of the algebraic equations from the following Krylov subspace: 
\begin{align}
K_{\mathsf{A}^2; n}:=\mathrm{span}\{\bs{b}, \mathsf{A}^2\bs{b}, \cdots, \mathsf{A}^{2n-2}\bs{b}\}, 
\label{eq:preconditioned_krylov_space}
\end{align}
if one uses the zero vector as the initial estimate. Notice that, in constructing the subspace $K_{\mathsf{A}^2; n}$, one needs $2n-2$ matrix-vector multiplications for the dense matrix $\mathsf{A}$. On the other hand, when solving the unpreconditioned system $\mathsf{A}\bs{y}=\bs{b}$ via GMRES, one obtains the following Krylov subspace: 
\begin{align}
K_{\mathsf{A}; 2n-1}:=\mathrm{span}\{\bs{b}, \mathsf{A}\bs{b}, \mathsf{A}^2\bs{b}, \cdots, \mathsf{A}^{2n-3}\bs{b}, \mathsf{A}^{2n-2}\bs{b}\}
\label{eq:unpreconditioned_krylov_space}
\end{align}
with the same amount of computations as for constructing $K_{\mathsf{A}^2; n}$. Comparing \eqref{eq:preconditioned_krylov_space} and \eqref{eq:unpreconditioned_krylov_space}, one may notice that $K_{\mathsf{A}^2; n}$ is a proper subset of $K_{\mathsf{A}; 2n-1}$. This means that, after the same amount of calculations, the GMRES for the unpreconditioned system explores the solution from a broader subspace than that for the preconditioned one. In other words, the GMRES for the unpreconditioned system always converges faster than that for the preconditioned one for a given tolerance.

In summary, we just need to appropriately arrange the integral operator $\cal A$ such that it satisfies ${\cal A}^2=\cal I + \cal K$ and to solve the \textcolor{black}{collocation-}discretised linear equations without multiplying the preconditioner. The present paper is based on the discovery that we can construct such a $\cal A$ for the Burton--Miller-based integral equation even for transmission problems. 

\subsection{Spectral property of the conventional Burton--Miller method for transmission problems}
Before deriving well-conditioned boundary integral equations based on the idea in the preceding section, let us review the conventional Burton--Miller-type boundary integral equation for the transmission problem in this section. Let us also share that the conventional equation is ill-conditioned.

To make the discussion simple, let us here consider the case of $M=2$, i.e. we have a single scatterer filled in a bounded domain $\Omega_2$ in the host matrix $\Omega_1$. Also assume that $\partial \Omega_2$ is simply connected and has no corners. In this case, a natural integral formulation based on the Burton--Miller method would couple the following integral equations:
\begin{align}
\frac{1}{2}u(\bs{x})+\frac{\alpha}{2} q(\bs{x})
&=
u^\mathrm{in}(\bs{x}) + \alpha q^\mathrm{in}(\bs{x})
+[{\cal S}_1q](\bs{x})
+\alpha [{\cal D}_1^*q](\bs{x})
-[{\cal D}_1 u](\bs{x})
-\alpha [{\cal N}_1 u](\bs{x}), 
\label{eq:bm_orig1} \\
\frac{1}{2}u(\bs{x})
&=
-[{\cal S}_2 q](\bs{x})
+[{\cal D}_2 u](\bs{x}), 
\label{eq:bm_orig2}
\end{align}
for $\bs{x}\in \Gamma:=\partial \Omega_1 \cup \partial \Omega_2$, where $q:=\nabla u\cdot \bs{n}$ and $q^\mathrm{in}:=\nabla u^\mathrm{in}\cdot \bs{n}$ are defined. Also, the symbols in calligraphic style indicate the following integral operators: 
\begin{align}
 [{\cal S}_i q](\bs{x}) := \int_{\partial \Omega_i} G_i(\bs{x}-\bs{y}) q(\bs{y}) \, \mathrm{d}\Gamma (\bs{y}), \label{eq:stmp} \\
 [{\cal D}_i u](\bs{x}) := \int_{\partial \Omega_i} \frac{\partial G_i(\bs{x}-\bs{y})}{\partial \bs{n}(\bs{y})} u(\bs{y}) \, \mathrm{d}\Gamma (\bs{y}), \label{eq:dtmp} \\
 [{\cal D}_i^* q](\bs{x}) := \int_{\partial \Omega_i} \frac{\partial G_i(\bs{x}-\bs{y})}{\partial \bs{n}(\bs{x})} q(\bs{y}) \, \mathrm{d}\Gamma (\bs{y}), \label{eq:dstmp} \\
 [{\cal N}_i u](\bs{x}) := \mathrm{p.\,f.}\,\int_{\partial \Omega_i} \frac{\partial^2 G_i(\bs{x}-\bs{y}) }{\partial \bs{n}(\bs{x}) \partial \bs{n}(\bs{y})} u(\bs{y}) \, \mathrm{d}\Gamma (\bs{y}), \label{eq:ddtmp}
\end{align}
where $\partial/\partial \bs{n}(\bs{x}):=\bs{n}(\bs{x})\cdot \nabla$ is the normal derivative, and ``p.\,f.'' indicates the finite part of the diverging integral. $G_i$ of the kernels in \eqref{eq:stmp}--\eqref{eq:ddtmp} is the fundamental solution of the two-dimensional Helmholtz equation and given as
\begin{align}
 G_i(\bs{x})=\frac{\mathrm{i}}{4} H_0^{(1)}(k_i|\bs{x}|), 
 \label{eq:fundamental_sol}
\end{align}
where $H_0^{(1)}$ is the Hankel function of the first kind and the $0^\mathrm{th}$-order. The integral equations in \eqref{eq:bm_orig1} and \eqref{eq:bm_orig2} are respectively obtained by taking the trace of the integral representations for $\bs{x}\in\Omega_1$ and  $\bs{x}\in\Omega_2$, respectively. For the former one, we employed the Burton--Miller method to make the fictitious eigenfrequency as far away from the real axis as possible. Note that the magnitude of the imaginary part of the eigenvalue can roughly be controlled by appropriately setting the coefficient $\alpha$ in \eqref{eq:bm_orig1}. We here use a widely accepted setting as $\alpha=-\mathrm{i}/k_1$~\cite{zheng2015burton}. 

By referring to the boundary conditions \eqref{eq:jump_in_u} and \eqref{eq:jump_in_w}, one may naturally obtain the following system of the boundary integral equations: 
\begin{align}
\begin{pmatrix}
\frac{1}{2}{\cal I}+{\cal D}_1+\alpha {\cal N}_1 & \alpha\left(\frac{1}{2}{\cal I}-{\cal D}^*_1\right)-{\cal S}_1 \\
\frac{1}{2}{\cal I}-{\cal D}_2 & \varepsilon_2 {\cal S}_2
\end{pmatrix}
\begin{pmatrix}
u\\
w
\end{pmatrix}
=
\begin{pmatrix}
u^\mathrm{in}+\alpha q^\mathrm{in}\\
0
\end{pmatrix}, 
\label{eq:naive_bie}
\end{align}
where ${\cal I}$ is the identity operator. The integral operator in \eqref{eq:naive_bie} is notoriously ill-conditioned. The upper and the lower diagonal parts are dominated by the hypersingular operator ${\cal N}$ and the \textcolor{black}{weakly singular} one ${\cal S}$, respectively. The eigenvalues of the former grow in magnitude, while those of the latter accumulate at the origin in the complex plane. This high contrast in the eigenvalues makes the resulting coefficient matrix considerably ill-conditioned. In \textcolor{black}{Section \ref{sec:2dom}}, we modify the integral equation \eqref{eq:naive_bie} with the help of the Calderon formulae to have a better conditioned one.

\section{\textcolor{black}{Calderon's formulae and their collocation discretisation}}
Before deriving an well-conditioned Burton--Miller equations for transmission problems, let us review the properties of the integral operators \eqref{eq:stmp}--\eqref{eq:ddtmp}. \textcolor{black}{First, the operators ${\cal S}_i$, ${\cal D}_i$, and ${\cal N}_i$ are continuous mappings for $H^{-1/2}(\partial \Omega_i)\rightarrow H^{1/2}(\partial \Omega_i)$, $H^{1/2}(\partial \Omega_i)\rightarrow H^{1/2}(\partial \Omega_i)$, $H^{1/2}(\partial \Omega_i)\rightarrow H^{-1/2}(\partial \Omega_i)$, respectively, and ${\cal D}_i^*$ is the adjoint of ${\cal D}_i$~\cite{coltonkress} for $i=1,\cdots, M$. If the boundary is smooth, the single- and double-layer operators are compact from $H^{1/2}(\partial\Omega_i)$ into $H^{1/2}(\partial \Omega_i)$~\cite{coltonkress}, and ${\cal N}_i$ is bounded.} It is also well known that the operators have the following relations: 
\begin{align}
\frac{\cal I}{4}&=-{\cal S}_i{\cal N}_i+{\cal D}_i^2 \label{eq:calderon1}, \\
0 &= {\cal S}_i {\cal D}^*_i - {\cal D}_i {\cal S}_i \label{eq:calderon2}, \\
0 &= -{\cal D}^*_i {\cal N}_i  + {\cal N}_i {\cal D}_i \label{eq:calderon3}, \\
\frac{\cal I}{4} &= ({\cal D}^*_i)^2-{\cal N}_i{\cal S}_i \label{eq:calderon4}.
\end{align}
The relations \eqref{eq:calderon1}--\eqref{eq:calderon4} are known as the Calderon formulae. It is also known that the formulae are valid modulo a compact operator even when the wave numbers of two operators in a single term are different from each other~\textcolor{black}{\cite{antoine2008integral}}.

\textcolor{black}{Since the product of compact operators is also compact~\cite{coltonkress}, from \eqref{eq:calderon1} and \eqref{eq:calderon4}, one notices that ${\cal S}_i$ (resp. ${\cal N}_i$) is the pseudo inverse of ${\cal N}_i$ (resp. ${\cal S}_i$) up to the factor of $-1/4$ modulo a compact operator. This is only true for the case that $\partial \Omega_i$ is smooth. However, for the collocation-discretised matrices $\mathsf{S}_i$ and $\mathsf{N}_i$ with constant elements respectively from the operators ${\cal S}_i$ and ${\cal N}_i$, the matrix products $\mathsf{S}_i\mathsf{N}_i$ and $\mathsf{N}_i\mathsf{S}_i$ approximates the identity up to the constant factor even when the boundary includes some corner points. This is because, in the collocation discretisation, the collocation points never get on any corner point; the collocation point is usually set on the centre of a boundary element. To numerically elaborate this fact, we computed the boundary element matrices $\mathsf{S}$ and $\mathsf{N}$ obtained as the collocation discretised version of the above operators on a square of unit side length. In the computation, each side of the square is divided into 100 line segments of equal length, the collocation points are set on the centre of boundary elements, and the wave number is set as $k=1.5$. The magnitude of element-wise absolute error of $-\mathsf{SN}$ and $-\mathsf{NS}$ to the identity matrix scaled by $1/4$ are at most $5.394\times 10^{-2}$ and $6.213\times 10^{-2}$, respectively. Based on this observation, it is concluded that the collocation boundary element matrices naturally inherit the spectral properties of integral operators.}

\section{The present formulation: A single material case}\label{sec:2dom}
We shall modify the original Burton--Miller integral equations \eqref{eq:naive_bie} so that the square of the operator $\cal A$ on the left-hand side does not involve the hypersingular operator $\cal N$ by exploiting the Calderon formulae. In other words, we shall use the Calderon formulae to eliminate the hypersingular operator $\cal N$ from the expression of ${\cal A}^2$. To this end, we first simply swap the first and second rows of the original equation \eqref{eq:naive_bie} as
\begin{align}
\begin{pmatrix}
\frac{1}{2}{\cal I}-{\cal D}_2 & \varepsilon_2 {\cal S}_2\\
\frac{1}{2}{\cal I}+{\cal D}_1+\alpha {\cal N}_1 & \alpha\left(\frac{1}{2}{\cal I}-{\cal D}^*_1\right)-{\cal S}_1
\end{pmatrix}
\begin{pmatrix}
u\\
w
\end{pmatrix}
=
\begin{pmatrix}
0\\
u^\mathrm{in}+\alpha q^\mathrm{in}
\end{pmatrix},
\label{eq:calderon_bie1}
\end{align}
whose integral operator in the left-hand side shall henceforth be denoted as ${\cal A}_1$. To see the ``preconditioned'' integral equation \eqref{eq:calderon_bie1} is well-conditioned, let us evaluate ${\cal A}_1^2$ as
\begin{align}
{\cal A}_1^2
&= 
\begin{pmatrix}
\frac{1}{4}{\cal I}+\alpha\varepsilon_2{\cal S}_2{\cal N}_1
&
0
\\
\frac{1+\alpha}{4}{\cal I}+\frac{\alpha(1+\alpha)}{2}{\cal N}_1-\alpha({\cal N}_1{\cal D}_2+\alpha {\cal D}_1^*{\cal N}_1)-\alpha {\cal S}_1{\cal N}_1
&
\frac{\alpha^2}{4}{\cal I}+\alpha\varepsilon_2 {\cal N}_1{\cal S}_2
\end{pmatrix}
+ {\cal K} \notag \\
&=
\begin{pmatrix}
\frac{1-\alpha\varepsilon_2}{4}{\cal I}& 0 \\
\frac{1+2\alpha}{4}{\cal I}+\frac{\alpha(1+\alpha)}{2}{\cal N}_1-\alpha({\cal N}_1{\cal D}_2+\alpha {\cal D}_1^*{\cal N}_1)
&
\frac{\alpha(\alpha-\varepsilon_2)}{4}{\cal I}
\end{pmatrix}
+ {\cal K},
\label{eq:A_1^2}
\end{align}
where we abusively used ${\cal K}$ to denote two compact operators different from each other. The last equality of \eqref{eq:A_1^2} is due to the Calderon formulae \eqref{eq:calderon1}--\eqref{eq:calderon4}. Although ${\cal A}_1^2$ still involves the hypersingular operator ${\cal N}_1$, owing to the block lower structure (except for a compact perturbation), its eigenvalue accumulates only at the following two points in the complex plane: $\frac{1-\alpha\varepsilon_2}{4}$ and $\frac{\alpha(\alpha-\varepsilon_2)}{4}$, and they are away from the origin of the complex plane (Recall that the setting $\alpha=-\mathrm{i}/k_1$ is used here). With these observations, we may expect that the GMRES applied to the linear algebraic equations obtained by discretising \eqref{eq:calderon_bie1} \textcolor{black}{via collocation} converges fast. Nonetheless, we cannot deny that the hypersingular operator ${\cal N}_1$ in ${\cal A}_1^2$ may affect the convergence of the linear solver. 

To derive a better conditioned integral equation \textcolor{black}{than \eqref{eq:calderon_bie1}}, we multiply a complex constant to \textcolor{black}{its} first row as
\begin{align}
\begin{pmatrix}
\beta\left(\frac{1}{2}{\cal I}-{\cal D}_2\right) & \beta \varepsilon_2 {\cal S}_2\\
\frac{1}{2}{\cal I}+{\cal D}_1+\alpha {\cal N}_1 & \alpha\left(\frac{1}{2}{\cal I}-{\cal D}^*_1\right)-{\cal S}_1
\end{pmatrix}
\begin{pmatrix}
u\\
w
\end{pmatrix}
=
\begin{pmatrix}
0\\
u^\mathrm{in}+\alpha q^\mathrm{in}
\end{pmatrix}, 
\label{eq:calderon_bieb}
\end{align}
with $\beta\in\mathbb{C}$. We henceforth denote the integral operator in \eqref{eq:calderon_bieb} as ${\cal A}_\beta$. Note that this notation is consistent with the previous one for the operator in \eqref{eq:calderon_bie1}; ${\cal A}_1$ is obtained by substituting $\beta=1$ to ${\cal A}_\beta$. As before, let us evaluate the square of the operator, yielding 
\begin{align}
{\cal A}_\beta^2= 
\begin{pmatrix}
\frac{\beta(\beta-\alpha\varepsilon_2)}{4}{\cal I} & 0 \\
\frac{\beta+2\alpha}{4}{\cal I}+\frac{\alpha(\beta+\alpha)}{2}{\cal N}_1-\alpha(\beta {\cal N}_1{\cal D}_2 + \alpha {\cal D}_1^*{\cal N}_1)
&
\frac{\alpha(\alpha-\beta\varepsilon_2)}{4}{\cal I}
\end{pmatrix}
+ {\cal K}, 
\label{eq:A_beta^2}
\end{align}
where ${\cal K}$ is again a compact operator. One may notice that setting $\beta=-\alpha$ in \eqref{eq:A_beta^2} eliminates all the hypersingular operators in the operator as follows: 
\begin{align}
{\cal A}_{-\alpha}^2=
\begin{pmatrix}
\frac{\alpha^2(1+\varepsilon_2)}{4} {\cal I} & 0\\
\frac{\alpha}{4}{\cal I}
& \frac{\alpha^2(1+\varepsilon_2)}{4}{\cal I}
\end{pmatrix}
+ {\cal K}, 
\label{eq:A_-alpha^2}
\end{align}
owing to the Calderon formula~\eqref{eq:calderon3}. Note also that ${\cal A}_{-\alpha}^2$ has the only one eigenvalue accumulation point at $\frac{\alpha^2(1+\varepsilon_2)}{4}$, while ${\cal A}_1^2$ has two. This beautiful property further improves the condition of the integral equation. Note that the clustering point may fall into the origin of the complex plane in the case of $\varepsilon_2=-1$. Such a situation may occur for a specially designed artificial material (e.g. metamaterial). The material constant is, however, often positive for natural materials.

\section{The present formulation: Multi-material cases}\label{sec:multi-material_cases}
\begin{figure}[h]
 \centering
 \includegraphics[scale=0.7]{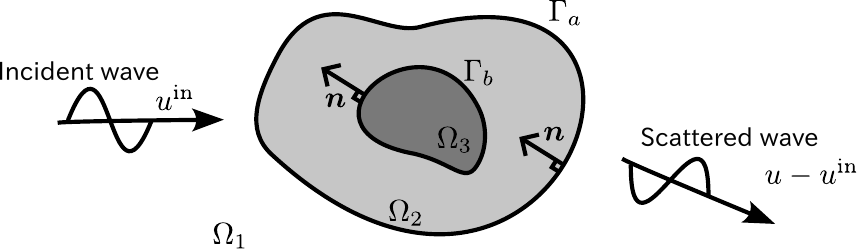}
 \caption{Illustrative sketch for multi-material scattering. The scatterer is composed of $\Omega_2$ and $\Omega_3$ each of which is filled with a different material.} \label{fig:3dom}
\end{figure}

In some applications, e.g. locally resonant sonic materials~\cite{liu2000locally}, it is of engineering interest to analyse wave scattering from composite materials. We are thus motivated to enhance our integral formulation~\eqref{eq:calderon_bieb} to solve such a problem. Here, we present our integral formulation for the setting illustrated in Figure \ref{fig:3dom}. Our scatterer is composed of two materials, each of which is filled in $\Omega_2$ and $\Omega_3$. Here, we assume that $\partial \Omega_3$ does not touch the host matrix $\Omega_1$. This assumption is just to simplify the presentation, but our formulation also works well for the case otherwise, as shall be seen in numerical experiments in Section \textcolor{black}{7}. It is also straightforward to apply the proposed formulation to scattering problems involving scatterers consisting of more than two materials. We introduce the following notation for boundaries in Figure \ref{fig:3dom}: $\Gamma_a:=\partial \Omega_1 \cap \partial \Omega_2$ and $\Gamma_b:=\partial \Omega_2 \cap \partial \Omega_3$, and denote the density functions living on the boundary $\Gamma_{a,b}$ as, for example, as $u_{a,b}$. We also set the unit normal on the boundaries directing into $\Omega_2$. 

With these settings, we use the Burton--Miller type boundary integral equation for the exterior domain $\Omega_1$ and the standard one for the other domains $\Omega_2$ and $\Omega_3$ as 
\begin{align}
\begin{cases}
\frac{1}{2}u_a(\bs{x})+\frac{\alpha}{2} w_a(\bs{x})
=
u^\mathrm{in}_a(\bs{x}) + \alpha q^\mathrm{in}_a(\bs{x})
+[{\cal S}_1w_a](\bs{x})
+\alpha [{\cal D}_1^*w_a](\bs{x})
-[{\cal D}_1 u_a](\bs{x})
-\alpha [{\cal N}_1 u_a](\bs{x}) & \bs{x}\in\Gamma_a \\
\frac{1}{2}u_a(\bs{x})
=
-\varepsilon_2[{\cal S}_2 w_a](\bs{x})-\varepsilon_2[{\cal S}_2 w_b](\bs{x})
+[{\cal D}_2 u_a](\bs{x})+[{\cal D}_2 u_b](\bs{x}) & \bs{x}\in\Gamma_a  \\
\frac{1}{2}u_b(\bs{x})
=
-\varepsilon_2[{\cal S}_2 w_a](\bs{x})-\varepsilon_2[{\cal S}_2 w_b](\bs{x})
+[{\cal D}_2 u_a](\bs{x})+[{\cal D}_2 u_b](\bs{x}) & \bs{x}\in\Gamma_b  \\
\frac{1}{2}u_b(\bs{x})
=
\varepsilon_3[{\cal S}_3 w_b](\bs{x})
-[{\cal D}_3 u_b](\bs{x}) & \bs{x}\in\Gamma_b
\end{cases}.\label{eq:3dom_bies}
\end{align}
We then need to appropriately couple the equations \eqref{eq:3dom_bies}. The conventional formulation may arrange the integral equations \eqref{eq:3dom_bies} as is, yielding the following equations:
\begin{align}
\begin{pmatrix}
\frac{1}{2}{\cal I}+{\cal D}_1+\alpha {\cal N}_1 &
0 &
\alpha\left(\frac{1}{2}{\cal I}-{\cal D}_1^*\right)-{\cal S}_1 &
0 \\
\frac{1}{2}{\cal I}-{\cal D}_2 &
-{\cal D}_2 &
\varepsilon_2 {\cal S}_2 & 
\varepsilon_2 {\cal S}_2 \\
-{\cal D}_2 &
\frac{1}{2}{\cal I}-{\cal D}_2 &
\varepsilon_2 {\cal S}_2 & 
\varepsilon_2 {\cal S}_2\\
0 &
\frac{1}{2}{\cal I}+{\cal D}_3 &
0 &
-\varepsilon_3 {\cal S}_3\\
\end{pmatrix}
\begin{pmatrix}
u_a\\
u_b\\
w_a\\
w_b\\
\end{pmatrix}
=
\begin{pmatrix}
u^\mathrm{in}_a+\alpha q^\mathrm{in}_a\\
0\\
0\\
0
\end{pmatrix}.
\label{eq:3dom_conv_bie}
\end{align}
Similarly to \eqref{eq:naive_bie}, it is easy to see that \eqref{eq:3dom_conv_bie} is highly ill-conditioned. If we naively extend the integral formalism \eqref{eq:calderon_bieb} to the current problem, we multiply the complex constant $\beta=-\alpha$ to the last three equations in \eqref{eq:3dom_bies} and ``appropriately'' arrange the equations, which gives the following: 
\begin{align}
\begin{pmatrix}
\beta\left(\frac{1}{2}{\cal I}-{\cal D}_2\right) &
-\beta {\cal D}_2 &
\beta \varepsilon_2 {\cal S}_2 & 
\beta \varepsilon_2 {\cal S}_2 \\
-\beta {\cal D}_2 &
\beta\left(\frac{1}{2}{\cal I}-{\cal D}_2\right) &
\beta \varepsilon_2 {\cal S}_2 & 
\beta \varepsilon_2 {\cal S}_2 \\
\frac{1}{2}{\cal I}+{\cal D}_1+\alpha {\cal N}_1 &
0 &
\alpha\left(\frac{1}{2}{\cal I}-{\cal D}_1^*\right)-{\cal S}_1 &
0 \\
0 &
\beta \left(\frac{1}{2}{\cal I}+{\cal D}_3\right) &
0 &
-\beta \varepsilon_3 {\cal S}_3
\end{pmatrix}
\begin{pmatrix}
u_a\\
u_b\\
w_a\\
w_b\\
\end{pmatrix}
=
\begin{pmatrix}
0\\
0\\
u^\mathrm{in}_a+\alpha q^\mathrm{in}_a\\
0
\end{pmatrix}.
\label{eq:3dom_calderon_bie0}
\end{align}
Let us denote the integral operator in the left-hand side of \eqref{eq:3dom_calderon_bie0} as ${\cal A}_{\beta; \mathrm{orig}}$. Its square reads
\begin{align}
{\cal A}_{\beta;\mathrm{orig}}^2=\begin{pmatrix}
\frac{\beta(\beta-\alpha\varepsilon_2)}{4}{\cal I} &
0 &
0 &
0 \\
\textcolor{black}{0} &
\frac{\beta^2}{4}{\cal I} &
0 &
0 \\
\frac{\beta+2\alpha}{4}{\cal I} + \frac{\alpha(\beta+\alpha)}{2}{\cal N}_1-\alpha(\beta {\cal N}_1{\cal D}_2+\alpha {\cal D}_1^*{\cal N}_1) & 
\textcolor{black}{0}&
\frac{\alpha(\alpha-\beta\varepsilon_2)}{4}{\cal I}& 
\textcolor{black}{0}&\\
0 &
\frac{\beta^2}{4}{\cal I} &
0 &
0 
\end{pmatrix} + {\cal K}, 
\label{eq:A^2_beta_orig}
\end{align}
with a compact operator ${\cal K}$. Unfortunately, the right-bottom block of \eqref{eq:A^2_beta_orig} ends up with a zero operator regardless of the choice of $\beta$. Due to the zero block in the diagonal \textcolor{black}{and the block diagonal structure}, the eigenvalue of ${\cal A}^2 _{\beta; \mathrm{orig}}$ may cluster at zero in the complex plane, yielding a possibly large condition number after the discretisation \textcolor{black}{via collocation}. 

To improve the spectral property of \eqref{eq:3dom_calderon_bie0}, i.e., to make the eigenvalue accumulation points away from the origin of the complex plane, let us look back at its single-material counterpart \eqref{eq:calderon_bieb}. The biggest difference between the integral operators in \eqref{eq:calderon_bieb} and \eqref{eq:3dom_calderon_bie0} is that the multi-material version \eqref{eq:3dom_calderon_bie0} includes the standard boundary integral equation in the lower half. One may thus come up with, to make the multi-material equation resemble the single-material one, replacing the equation in the last row of \eqref{eq:3dom_calderon_bie0} by the one of the Burton--Miller type as
\begin{align}
\begin{pmatrix}
\beta\left(\frac{1}{2}{\cal I}-{\cal D}_2\right) &
-\beta {\cal D}_2 &
\beta \varepsilon_2 {\cal S}_2 &
\beta \varepsilon_2 {\cal S}_2 \\
-\beta {\cal D}_2 &
\beta\left(\frac{1}{2}{\cal I}-{\cal D}_2\right) &
\beta \varepsilon_2 {\cal S}_2 &
\beta \varepsilon_2 {\cal S}_2 \\
\frac{1}{2}{\cal I}+{\cal D}_1+\alpha {\cal N}_1 &
0 &
\alpha\left(\frac{1}{2}{\cal I}-{\cal D}_1^*\right)-{\cal S}_1 &
0 \\
0 &
\frac{1}{2}{\cal I}+{\cal D}_3+\gamma {\cal N}_3 &
0 &
\gamma\varepsilon_3\left(\frac{1}{2}{\cal I}-{\cal D}_3^*\right)-\varepsilon_3{\cal S}_3 &
\end{pmatrix}
\begin{pmatrix}
u_a\\
u_b\\
w_a\\
w_b\\
\end{pmatrix}
=
\begin{pmatrix}
0\\
0\\
u^\mathrm{in}_a+\alpha q^\mathrm{in}_a\\
0
\end{pmatrix},
\label{eq:3dom_calderon_bie}
\end{align}
where $\gamma\in\mathbb{C}$ is the coefficient for the newly introduced Burton--Miller equation. We henceforth denote the operator in \eqref{eq:3dom_calderon_bie} as ${\cal A}_{\beta,\gamma;\mathrm{mod}}$. We here determine the coefficients $\beta$ and $\gamma$ in ${\cal A}_{\beta,\gamma;\mathrm{mod}}$ so that ${\cal A}_{\beta,\gamma; \mathrm{mod}}^2$ involves the hypersingular operator as less as possible. This is achieved by setting the tunable parameters as $\beta=-\alpha$ and $\gamma=\frac{\alpha}{\varepsilon_3}$, which gives
\begin{align}
{\cal A}_{-\alpha,\alpha/\varepsilon_3;\mathrm{mod}}^2=
\begin{pmatrix}
\frac{\alpha^2(1+\varepsilon_2)}{4}{\cal I} &
\textcolor{black}{0}&
0 &
0 \\
\textcolor{black}{0}&
\frac{\alpha^2\left(1+\frac{\varepsilon_2}{\varepsilon_3}\right)}{4}{\cal I} &
0 &
0 \\
\frac{\alpha}{4}{\cal I} &
\textcolor{black}{0}&
\frac{\alpha^2(1+\varepsilon_2)}{4}{\cal I} &
\textcolor{black}{0}\\
\textcolor{black}{0}&
\frac{\alpha}{4}{\cal I} &
\textcolor{black}{0}&
\frac{\alpha^2\left(1+\frac{\varepsilon_2}{\varepsilon_3}\right)}{4}{\cal I}
\end{pmatrix} + {\cal K}.
\label{eq:A^2_beta_mod}
\end{align}
It is easy to see that the eigenvalue clusters are found only at \textcolor{black}{$\frac{\alpha^2(1+\varepsilon_2)}{4}$ and $\frac{\alpha^2(1+\varepsilon_2/\varepsilon_3)}{4}\neq 0$}. We may thus conclude that the proposed integral formalism \eqref{eq:3dom_calderon_bie} with $\beta=-\alpha$ and $\gamma=\alpha/\varepsilon_3$ is well-conditioned. Note that there is another possibility as the integral equations \eqref{eq:3dom_calderon_bie}; The equations in the second and fourth rows can be interchanged, \textcolor{black}{for example}. In that case, one needs to use the ordinary and the Burton--Miller formulations for the new second and fourth rows, respectively, to have a well-conditioned system. The parameters $\beta$ and $\gamma$ should also be adjusted so that \textcolor{black}{all} the hypersingular operators are eliminated from the square of the underlying operator. 

Although we have presented the formulation for a specific case shown in Figure \ref{fig:3dom}, applying the present formulation to more general settings is almost straightforward and can systematically be done. We here just show a recipe for constructing the integral equation corresponding to \eqref{eq:3dom_calderon_bie}: One first writes down all the necessary integral equations as \eqref{eq:3dom_bies} by taking the appropriate boundary traces of the integral representations for $u$. The equation from the infinite domain $\Omega_1$ should be of the Burton--Miller type. The boundary integral equations derived from the domain on whose boundary the normal is positively directed are also set as the Burton--Miller type. The coefficient for the Burton--Miller formulation thus set for $\Omega_i$ for $i\neq 1$ is set as $\alpha/\varepsilon_i$, where $\alpha=-\frac{i}{k_1}$ is that for $\Omega_1$. The rest of the integral equations are of standard type (i.e. consisting of the single- and double-layer potentials) multiplied by a complex constant $-\alpha$. The standard and the Burton--Miller type operators are arranged in the upper and the lower half of the integral operator, respectively.

\section{Numerical demonstrations}\label{sec:numerical_experiments}
In this section, we demonstrate the validity and efficiency of the proposed integral formalism through several numerical experiments. The boundary integral equations are discretised by the collocation method with constant elements. The boundary integrals are evaluated by the Gauss--Legendre quadrature with 10 Gaussian nodes. In the evaluation, the weak singularity for \eqref{eq:stmp} is handled with a standard ``subtraction and addition'' technique, and the integral with the hypersingular kernel \eqref{eq:ddtmp} is also regularised~\textcolor{black}{\cite{canino1998numerical}}. As the linear algebraic solver, the GMRES without restart is utilised, and its tolerance is set as $10^{-8}$. 

\subsection{Single material cases}
\begin{figure}[h]
 \centering
 \includegraphics[scale=0.5]{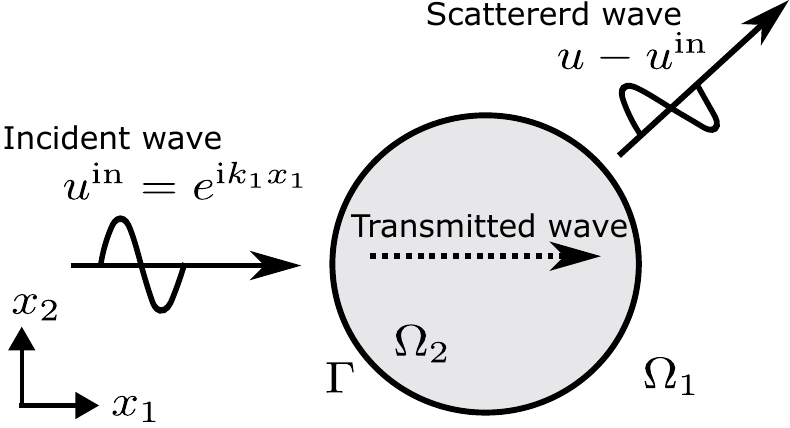}
 \caption{Benchmark setting for the single material case.}\label{fig:2dom_circle}
\end{figure}

Let us first consider the single material case. Here, as a benchmark, we use the origin-centred unit circular scatterer in Figure \ref{fig:2dom_circle}. The material constants for the host medium $\Omega_1$ and the inclusion $\Omega_2$ are set as $\varepsilon_1=1.0$ and $\varepsilon_2$, respectively. The domain is impinged by a plane incident wave $u^\mathrm{in}=e^{\mathrm{i}k_1 x_1}$ propagating in $x_1$-direction.

\begin{figure}[h]
 \centering
  \includegraphics[scale=0.6]{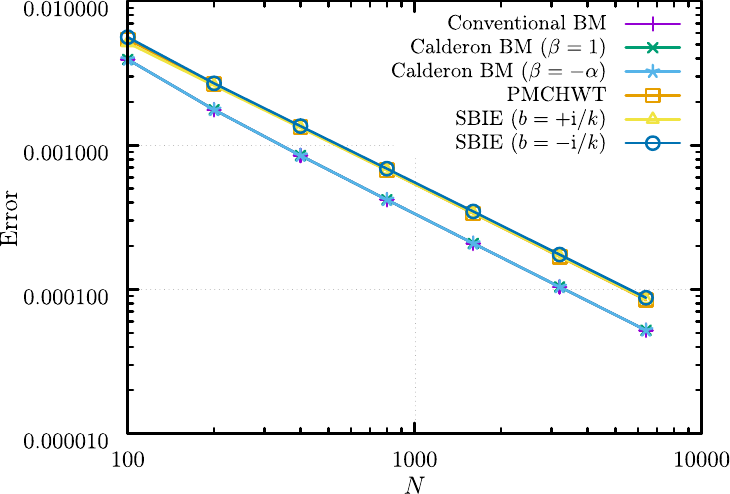}
  \caption{Relative $\ell_2$-error in $u$ vs the number of boundary elements $N$ for the case of single material scattering. ``Conventional BM'', ``Calderon BM ($\beta=1$)'', and ``Calderon BM ($\beta=-\alpha$)'' respectively solve the integral equations \eqref{eq:naive_bie}, \eqref{eq:calderon_bie1}, and \eqref{eq:calderon_bieb} with $\beta=-\alpha$. The same keys shall be used in the figures to follow.}\label{fig:00_2dom_nelem_bin_00fig_2dom_nelem_error}
\end{figure}

We first check the accuracy. To this end, we computed the total field $u$ at the 10,201~($=101\times 101$) field points equally distributed in $\{\bs{x} \mid -2.99 \le x_i \le 2.99,~i=1,~2 \}$ and compared the numerical results with the analytical ones. In this computation, the unit circle is approximated by its inscribed polygon with $N$ edges, and the angular frequency and the material constant for the inclusion are set as $\omega=5.0$ and $\varepsilon_2=2.0$, respectively. The analytical solution in this case is given as 
\begin{equation}
 u(\bs{x})=
\begin{cases}
u^\mathrm{in}(\bs{x})+\displaystyle\sum_{n=-\infty}^{\infty}a_n H_n^{(1)}(k_1|\boldsymbol{x}|)e^{in\theta} & \boldsymbol{x}\in\Omega_1 \\
\displaystyle\sum_{n=-\infty}^{\infty}b_n J_n(k_2|\boldsymbol{x}|)e^{in\theta} & \boldsymbol{x}\in\Omega_2
\end{cases}, 
\label{eq:2dom_anal}
\end{equation}
where $J_n$ is the Bessel function of $n^\mathrm{th}$-order, and $a_n$ and $b_n$ are complex constants determined by the boundary conditions \eqref{eq:jump_in_u} and \eqref{eq:jump_in_w}. In the practical computation, the infinite series in \eqref{eq:2dom_anal} are truncated by $\pm 50$. Figure \ref{fig:00_2dom_nelem_bin_00fig_2dom_nelem_error} shows the relative $\ell_2$-error between the real parts of the numerical and analytical solutions defined as
\begin{align}
 \mathrm{Error}:=\sqrt{\frac{\sum_{j=1}^{10201} \Re\left[u^\mathrm{num}(\bs{x}_j)-u^\mathrm{ana}(\bs{x}_j)\right]^2}{\sum_{j=1}^{10201} \Re\left[u^\mathrm{ana}(\bs{x}_j)\right]^2}}, 
\end{align}
where $\bs{x}_j$ indicates the $j^\mathrm{th}$-field point, and $u^\mathrm{num}$ and $u^\mathrm{ana}$ are the numerical and analytical solutions, respectively. 
The figure shows the results obtained by the proposed integral formulation \eqref{eq:calderon_bieb} with $\beta=-\alpha$ and $\beta=1$ and the conventional one \eqref{eq:naive_bie} as well as those obtained by widely used ones, i.e. the PMCHWT and SBIE formulations. The integral operators for the PMCHWT equation are ``appropriately'' arranged such that the square of the underlying operator is equal to the identity modulo a compact operator~\cite{niino2012preconditioning}. The SBIE is constructed by using the following ansatz for $u$ in $\bs{x}\in\Omega_1$:
\begin{align}
u(\bs{x}) = u^\mathrm{in}(\bs{x}) + [{\cal S}_1\mu](\bs{x}) + b [{\cal D}_1 \mu] (\bs{x}), 
\label{eq:sbie}
\end{align}
where $\mu$ is the unknown density, which is determined by the boundary conditions \eqref{eq:jump_in_u} and \eqref{eq:jump_in_w} and the Green identity for $\bs{x}\in\Omega_2$. $b$ in \eqref{eq:sbie} is a complex parameter having a similar role as the coefficient $\alpha$ for the Burton--Miller method, controlling the position of the fictitious eigenvalues. Here, $b$ is set either $b=\mathrm{i}/k_1$ or $-\mathrm{i}/k_1$. From the figure, we may observe that the accuracy of all the method scales with $N^{-1}$, which is a typical behaviour of the solution obtained by the collocation BEM for two-dimensional problems with constant elements. It is also observed that the Burton--Miller method gives slightly more accurate solutions than the other formulations: PMCHWT and SBIE. Furthermore, the preconditioning would, of course, not affect the accuracy.

\begin{figure}[h]
 \centering
  \includegraphics[scale=0.6]{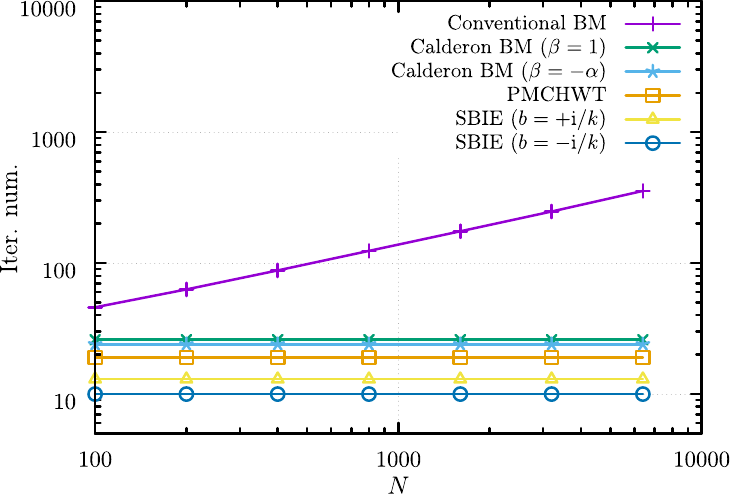}
  \caption{The number of iterations for GMRES vs $N$ for the case of single material scattering.}\label{fig:00_2dom_nelem_bin_01fig_2dom_nelem_iter.pdf}
\end{figure}

We then check the number of iterations of the GMRES against $N$ in Figure \ref{fig:00_2dom_nelem_bin_01fig_2dom_nelem_iter.pdf}. While the iteration number for ``Conventional BM'' grows up, that for the other formulations remains constant as $N$ increases. Note that the GMRES with $\beta=-\alpha$ always converges faster than that with $\beta=1$ as easily expected from the theoretical discussions in Section \ref{sec:2dom}; see \eqref{eq:A_1^2} and \eqref{eq:A_-alpha^2} specifically. It may, however, be surprising that the difference is not so significant. The hypersingular operators in the left-bottom block in \eqref{eq:A_1^2} might not have a serious effect on the efficiency; Recall that they would not affect the spectral property of the non-compact part of ${\cal A}_1^2$ owing to the zero operator in the upper-right block. The number of iterations for the present formulations is slightly larger than but almost comparable with that for the PMCHWT and SBIEs.

\begin{figure}[h]
 \centering
 \includegraphics[scale=0.55]{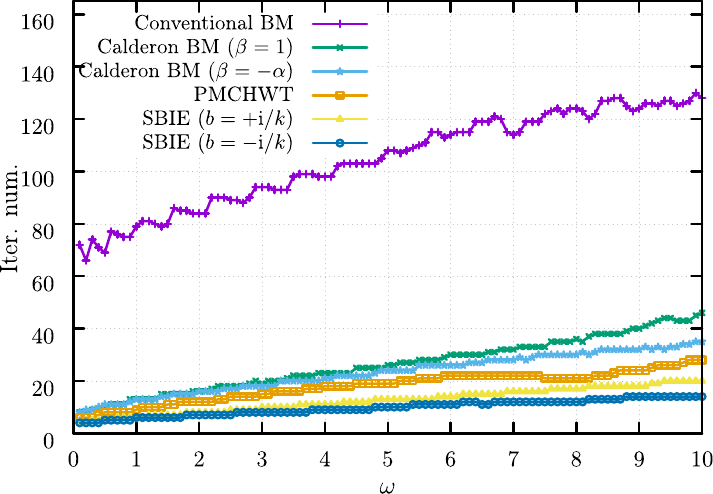}~~~~~~
 \includegraphics[scale=0.55]{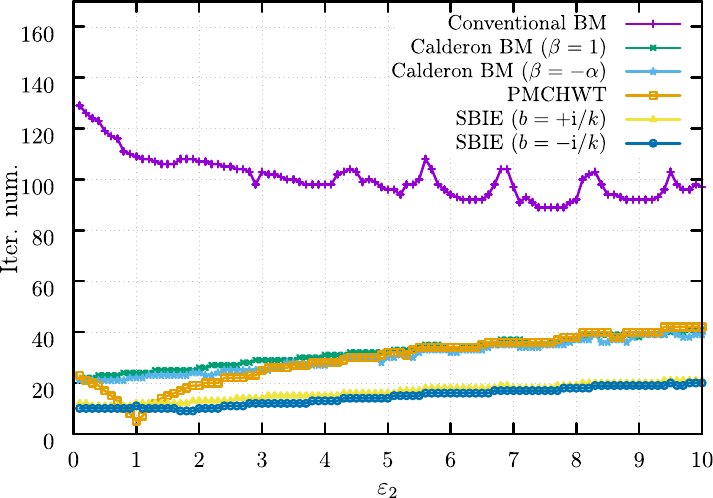}
 \caption{The number of GMRES iterations against the incident angular frequency $\omega$ (left) and the material constant of the circular inclusion $\varepsilon_2$ (right).}
 \label{fig:01_and_02}
\end{figure}

We then verify the efficiency dependence of our formulation on excitation angular frequency $\omega$ and the contrast ratio in material constants $\varepsilon_2~(= \varepsilon_2/\varepsilon_1)$ of the inclusion and the exterior medium. To this end, we again use the circular example (Figure \ref{fig:2dom_circle}) and measure the number of GMRES iterations for various $\omega$ and $\varepsilon_2$. The left and right of Figure \ref{fig:01_and_02} shows the number of iterations vs $\omega$ and $\varepsilon_2$, respectively. In the computation for the left figure, $\varepsilon_2$ is fixed as $2.0$. For right, $\omega$ is fixed as $5.0$. For both computations, the number of boundary elements $N$ is fixed as $N=600$. The same observations as the previous example are also seen in these examples; The proposed formulations as well as the PMCHWT and SBIEs outperform the conventional Burton--Miller method, and the proposed one with $\beta=-\alpha$ is almost always faster than the one with $\beta=1$, and the performance of the proposed methods slightly worse than the PMCHWT and SBIEs. The performance difference is, however, not significant. It is also worth noting that, in Figure \ref{fig:01_and_02}, while there are several peaks in the plot for the conventional Burton--Miller method, there are none for the rest of the methods. The peaks correspond to the resonance angular frequencies; At these frequencies, the transmitted wave is trapped in the inclusion $\Omega_2$. It is well-known that the linear system obtained from the boundary integral equation may considerably be ill-conditioned for such a case. Nonetheless, the proposed formulations (as well as the PMCHWT and SBIEs) robustly work even in such a case. With these observations, especially given the high accuracy (see Figure \ref{fig:00_2dom_nelem_bin_00fig_2dom_nelem_error}), the proposed method would also be attractive when solving transmission problems consisting of a single scatterer.

\begin{figure}[h]
 \centering
 \includegraphics[scale=0.6]{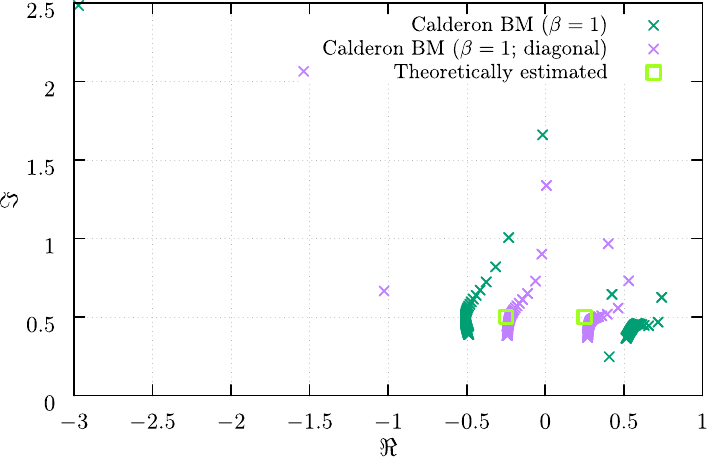}
 \caption{The eigenvalue distribution of $\mathsf{A}_1^2$ and its diagonal blocks.}\label{fig:eigenvalue_00_2dom_eigenvalue_w1_fig_2dom_eigenvalue3}
\end{figure}

To further corroborate that the fast GMRES convergence for the proposed methods originated from the carefully designed spectral properties of \eqref{eq:A_1^2} and \eqref{eq:A_-alpha^2}, we here check the eigenvalues of the discretised version of the operators, i.e. matrices $\mathsf{A}_1^2$ and $\mathsf{A}_{-\alpha}^2$. Here, the circular scatterer problem in Figure \ref{fig:2dom_circle} is again used, and its surface is discretised into $N=100$ boundary elements. Also, the angular frequency $\omega$ and the material constant $\varepsilon_2$ are set as $\omega=1.0$ and $\varepsilon_2=2.0$, respectively. With these settings, the eigenvalue accumulation point(s) of the operators ${\cal A}^2_1$ in \eqref{eq:A_1^2} and ${\cal A}^2_{-\alpha}$ in \eqref{eq:A_-alpha^2} are found at $\frac{1-\alpha\varepsilon_2}{4}, \frac{\alpha(\alpha-\varepsilon_2)}{4}=\pm\frac{1}{4}+\frac{1}{2}\mathrm{i}$ and $\frac{\alpha^2(1+\varepsilon_2)}{4}=-\frac{3}{4}+0\mathrm{i}$, respectively. 
\begin{figure}[h]
 \centering
 \includegraphics[scale=0.6]{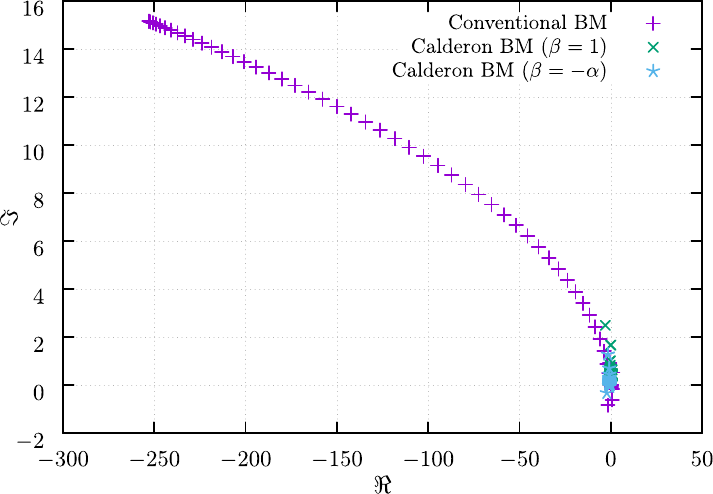}
 \caption{The eigenvalue distribution of (the discretised version of) the square of the operator in the conventional \eqref{eq:naive_bie}, the proposed integral equations \eqref{eq:calderon_bie1} and \eqref{eq:calderon_bieb} with $\beta=-\alpha$.}
 \label{fig:eigenvalue_00_2dom_eigenvalue_w1_fig_2dom_eigenvalue}
\end{figure}
Figure \ref{fig:eigenvalue_00_2dom_eigenvalue_w1_fig_2dom_eigenvalue} shows the eigenvalue distribution of $\mathsf{A}_1^2$ and $\mathsf{A}_{-\alpha}^2$. The distribution of the square of the conventional integral equation \eqref{eq:naive_bie} is also plotted in the same figure. The eigenvalues for the proposed methods are clustered in the vicinity of the origin, while those for the conventional one diverge. These spectral properties for the proposed operators contributed to reducing the number of GMRES iterations. To further explore the differences between the $\beta=1$ and $\beta=-\alpha$ cases, we show in Figure \ref{fig:eigenvalue_00_2dom_eigenvalue_w1_fig_2dom_eigenvalue2} a magnified view around the origin of Figure \ref{fig:eigenvalue_00_2dom_eigenvalue_w1_fig_2dom_eigenvalue}. Note that all eigenvalues for the relevant matrices are plotted in the figure. 
\begin{figure}[h]
 \centering
 \includegraphics[scale=0.6]{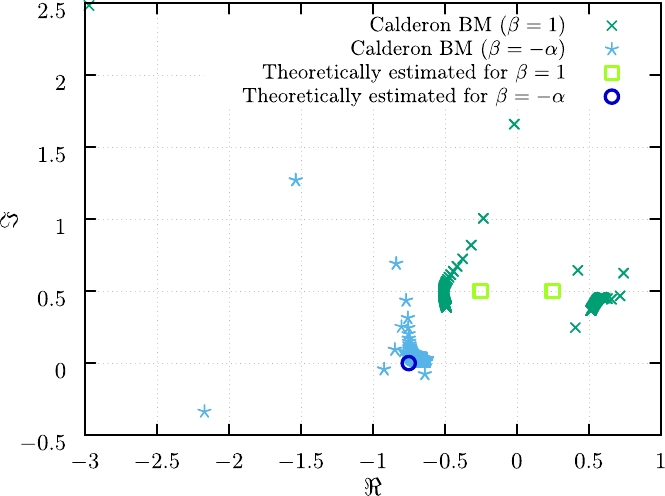}
 \caption{The eigenvalue distribution of $\mathsf{A}_1^2$ and $\mathsf{A}_{-\alpha}^2$ (The close-up of Figure \ref{fig:eigenvalue_00_2dom_eigenvalue_w1_fig_2dom_eigenvalue}).}
 \label{fig:eigenvalue_00_2dom_eigenvalue_w1_fig_2dom_eigenvalue2}
\end{figure}
As expected, the plots for $\beta=1$ cluster around two distinct locations, while those for $\beta=-\alpha$ converge toward a single point. This difference may explain the better performance of $\beta=-\alpha$ setting than $\beta=1$ one. One may also notice that, while the eigenvalues of ${\cal A}_{-\alpha}^2$ clusters near the theoretically estimated position, those of ${\cal A}_1^2$ do a bit different locations from the estimation. This inconsistency for ${\cal A}_1^2$ may caused by the hypersingular operators in \eqref{eq:A_1^2} in the left-bottom block. Because the upper-right block of ${\cal A}_1^2$ is not exactly zero but includes a compact perturbation, the hypersingular operators may have some impact on the eigenvalues of the discretised matrix. To check the impact, we here plot the eigenvalue distribution of the (discretised version of) diagonal blocks of ${\cal A}_1^2$ in Figure \ref{fig:eigenvalue_00_2dom_eigenvalue_w1_fig_2dom_eigenvalue3}. As expected, the eigenvalues of the diagonal blocks are gathered around the expected locations. Nonetheless, the plots for $(\beta=1)$ and $(\beta=1; \mathrm{diagonal})$ are close to each other, from which we may conclude that the effect of ${\cal N}$s in the left-bottom block is not so catastrophic. 

\subsection{Multi-material cases}
\begin{figure}[h]
 \centering
 \includegraphics[scale=0.5]{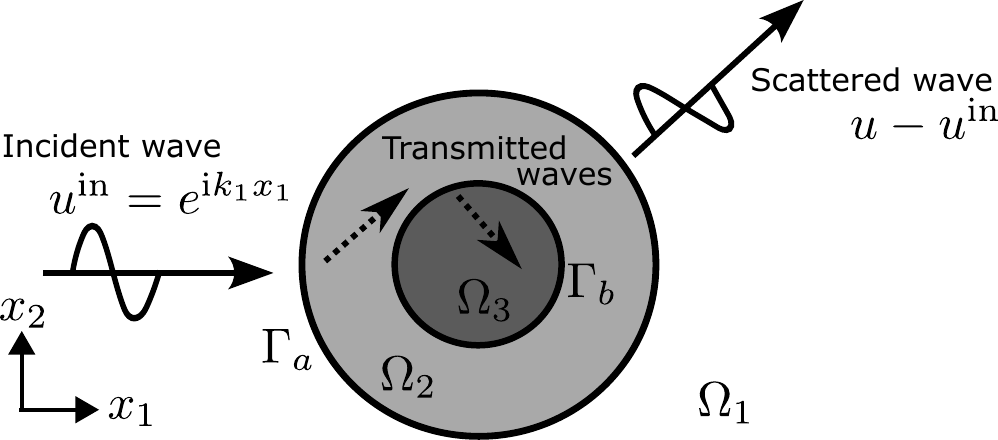}
 \caption{Illustrative sketch for the scattering by concentric circular materials.}\label{fig:3dom_circle}
\end{figure}

We then check the performance of the proposed methods for transmission problems defined in a domain with multiple materials. In the first example, a benchmark problem shown in Figure \ref{fig:3dom_circle} is used; In the exterior host medium $\Omega_1$ with $\varepsilon_1=1$, origin-centred concentric circular materials are allocated. The outer circle of radius 2.0 with $\varepsilon_2=2.0$ and inner one of radius 1.0 with $\varepsilon_3=3.0$ are denoted as $\Omega_2$ and $\Omega_3$, respectively. We again consider the plane wave incidence as $u^\mathrm{in}(\bs{x})=e^{\mathrm{i}k_1x_1}$ for $\bs{x}\in\Omega_1$. 

We solve this problem by the conventional Burton--Miller integral equation \eqref{eq:3dom_conv_bie}, the equation \eqref{eq:3dom_calderon_bie0} with $\beta=-\alpha$ obtained by naively extending the proposed formulation for single-material case, and the one tailored for multi-material case \eqref{eq:3dom_calderon_bie} with $\beta=-\alpha$ and $\gamma=\alpha/\varepsilon_3$ by utilising the Burton--Miller method for $\Omega_3$. The results by them are labelled as ``Conventional BM'', ``Calderon BM ($\beta=-\alpha$)'', and ``Calderon BM (mod.)'' in the figures to follow.
\begin{figure}[h]
 \centering
 \includegraphics[scale=0.55]{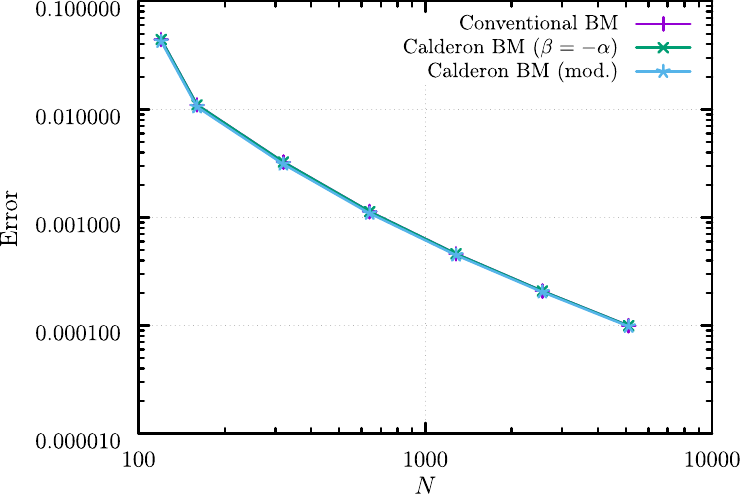}~~~~~~
 \includegraphics[scale=0.55]{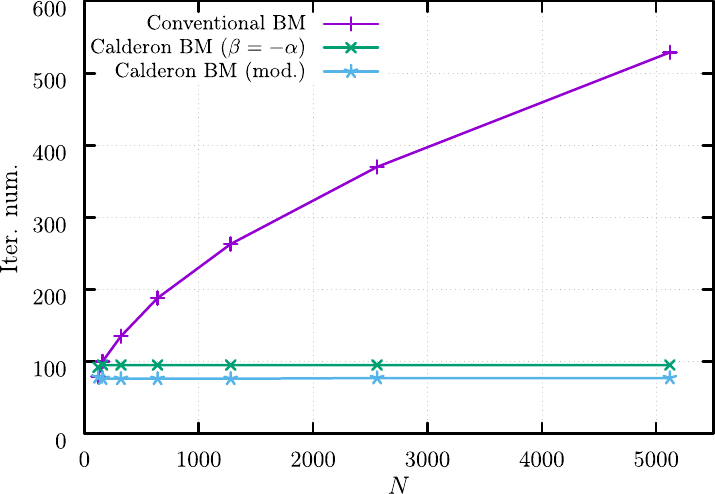}
 \caption{The accuracy (left) and the number of GMRES iterations (right) vs the number of boundary elements $N$ in the case of concentric circular scatterers.}\label{fig:04_3dom_nelem}
\end{figure}
In the computations, the inner and outer circles are respectively discretised $N/3$ and $2N/3$ constant boundary elements, and $u$ at 10,201 field points distributed in $\{ \bs{x} \mid -2.99 \le x_i \le 2.99,~i=1,2 \}$ of equal intervals are evaluated. The left and right figures of Figure \ref{fig:04_3dom_nelem} show the relative $\ell_2$-error of the numerical solutions to the analytical ones and the number of GMRES iterations versus $N$, respectively. The analytical solution for this case can also be derived similarly to the single-material case \eqref{eq:2dom_anal}. Observing Figure \ref{fig:04_3dom_nelem} (left), we found that the proposed preconditioning has no significant impact on accuracy. It is important to highlight that, for ``Calderon BM (mod.)'' unlike the cases of ``Calderon BM ($\beta=-\alpha$)'' and the single-material formulations, we do not only ``precondition'' the equation but also make a change to the integral formulation itself. Specifically, the integral equation derived from $\Omega_3$ is replaced by the Burton--Miller one. While such a modification could potentially affect the accuracy of the integral formulation, our investigation suggests that this concern is, at least in this example, not applicable. Figure \ref{fig:04_3dom_nelem} (right) shows that the GMRES iterations for the proposed Calderon-preconditioned systems remain constant for various $N$, while that of the conventional one does depend on the model size $N$. One also observes that the modified version \eqref{eq:3dom_calderon_bie} always outperforms the original one \eqref{eq:3dom_calderon_bie0}. 

\begin{figure}[h]
 \centering
 \includegraphics[scale=0.55]{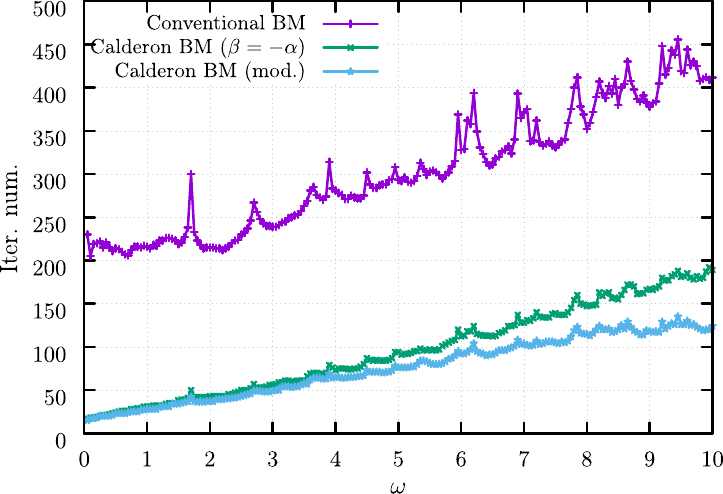}~~~~~~
 \includegraphics[scale=0.55]{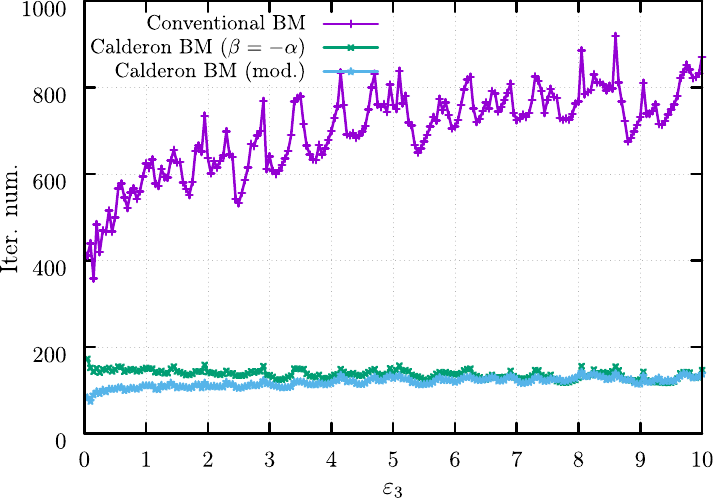}
 \caption{The number of GMRES iterations for the concentric circular scattering. The left and right show the dependence of the iteration numbers on the excitation angular frequency and the material constants, respectively.}\label{fig:3dom_omega_and_epsln}
\end{figure}

We then check the performance dependence of the proposed methods on the incident angular frequency $\omega$ and the material constants $\varepsilon_2$ and $\varepsilon_3$. For the former, we fix the material constants as before: $\varepsilon_i = i$~for $i=1,\cdots, 3$. For the latter, we fix $\omega=5.0$ and sweep the inclusion material constants as $\varepsilon_2=1/\varepsilon_3$ for $0<\varepsilon_3\le 10$. For both computations, the inner and outer circular surfaces are discretised into 600 and 1200 constant boundary elements, respectively. Figure \ref{fig:3dom_omega_and_epsln} shows the number of GMRES iterations for conventional and proposed BEMs of the Burton--Miller type.
According to the figures, the proposed methods can reduce the number of iterations. The sharp peaks for the conventional formulation are also considerably suppressed for the proposed formulations, which means that the Calderon-preconditioned integral equations work robustly against the physical resonance. Also, the modified version (mod.) almost always converges faster than the original ($\beta=-\alpha$) in these cases. When $\varepsilon_3$ is large (and correspondingly $\varepsilon_2$ is small, and their contrast ratio is high), the original and modified versions show a comparable performance. This might be explained by investigating the accumulation points of the eigenvalues for ${\cal A}^2_{-\alpha; \mathrm{orig}}$ in \eqref{eq:A^2_beta_orig} and ${\cal A}^2_{-\alpha,\alpha/\varepsilon_3; \mathrm{mod}}$ in \eqref{eq:A^2_beta_mod}. The former and latter accumulate, in the limit of $\varepsilon _3\rightarrow \infty$, towards the two points $-1/100$ and $0$ and the single point $-1/100$, all of which are in close proximity.

\begin{figure}[h]
 \centering
 \includegraphics[scale=0.6]{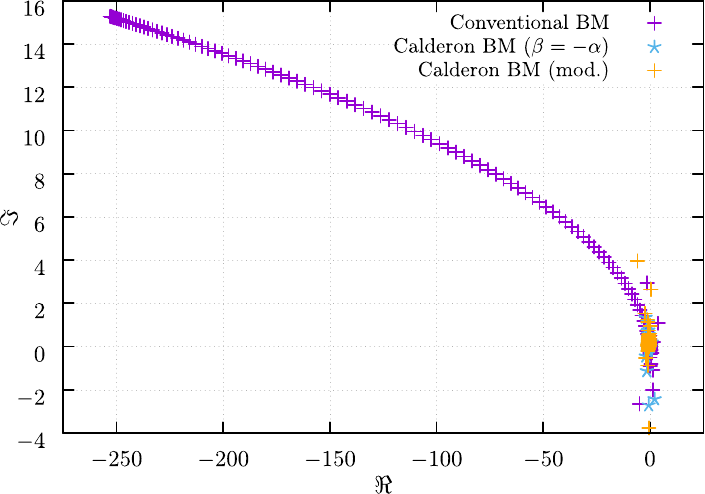}~~~~~~
 \includegraphics[scale=0.6]{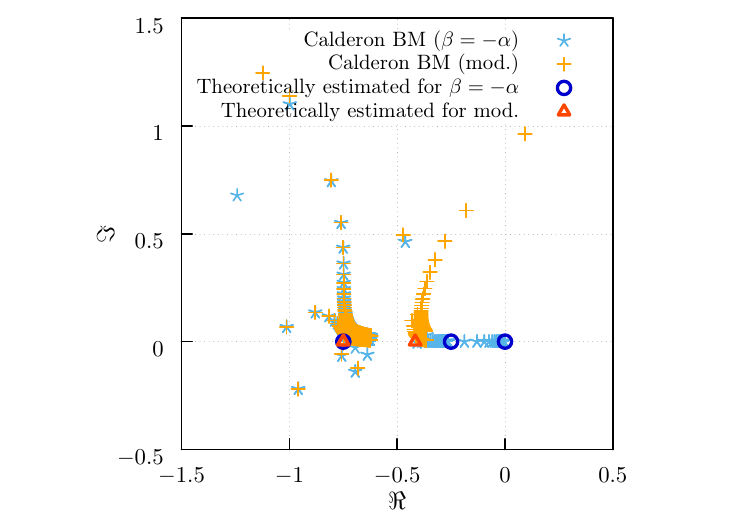}
 \caption{\textcolor{black}{The eigenvalue distributions of the square of (the discretised version of) the coefficient matrices of \eqref{eq:3dom_bies}, \eqref{eq:3dom_calderon_bie0} with $\beta=-\alpha$, and the modified one \eqref{eq:3dom_calderon_bie} with $\beta=-\alpha, \gamma=\alpha/\varepsilon_3$. The right figure shows the close-up near the origin of the left.}}
\label{fig:3dom_eigenvalue}
\end{figure}

We then check the eigenvalue distribution of the square of the coefficient matrix of each formulation. For this, we fixed the number of boundary elements on the inner and outer boundaries as 100 and 200, respectively. The material constants and the angular frequency are set as $(\varepsilon_1, \varepsilon_2, \varepsilon_3)=(1.0, 2.0, 3.0)$ and $\omega=1.0$, respectively. Figure \ref{fig:3dom_eigenvalue} shows the results.
The eigenvalues of the square of the proposed coefficint matrices cluster around the origin, while those of the conventional formulation diverge. These spectral properties account for the fast GMRES convergence of the proposed formulations. The right of Figure \ref{fig:3dom_eigenvalue} shows the close-up near the origin of the left one. One finds that the eigenvalues for the proposed formulations cluster at a few points. Note here that not all the eigenvalues for ($\beta=-\alpha$) and (mod.) are plotted in the right figure. The range is truncated so that one easily understands the behaviour of the eigenvalues in the vicinity of their clustering points. The eigenvalues with the largest magnitude for ($\beta=-\alpha$) and (mod.) are found at $2.1706-2.4538\mathrm{i}$ and $-5.9178+3.9583\mathrm{i}$, respectively. In the right figure (the magnified one), the theoretically estimated accumulation points for the integral operators are also plotted. \textcolor{black}{One confirms that the eigenvalues accumulate around the estimated points.} As expected, the original version of the Calderon-preconditioned equation with $\beta=-\alpha$ involves the coefficient matrix whose square has eigenvalues with a tiny magnitude. On the other hand, the coefficient matrix of the modified version does not. These spectral properties explain the performance difference of the proposed solvers, concluding that the modified version is superior to the original one. 

\begin{figure}[h]
 \centering
 \includegraphics[scale=0.5]{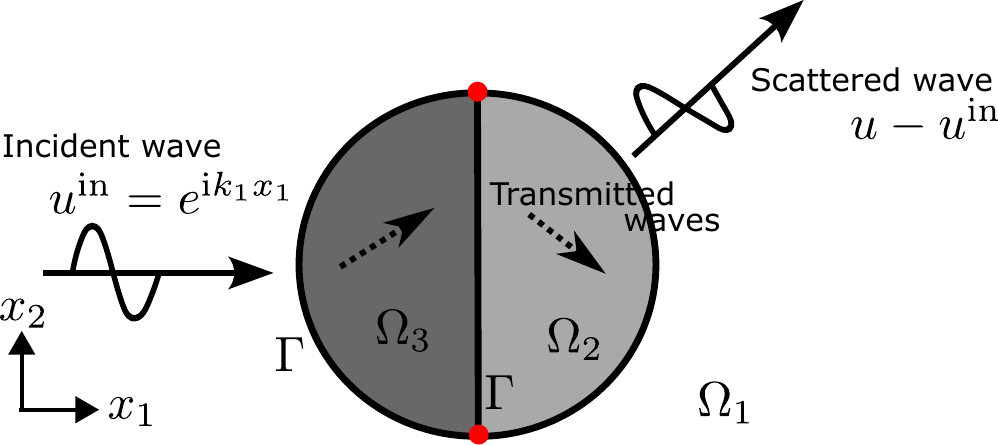}
 \caption{Transmission and scattering by a scatterer with material junctions, where three materials converge at a single point (indicated by red marks).}
 \label{fig:3mata_circle}
\end{figure}

\begin{figure}[h]
 \centering
 \includegraphics[scale=0.35]{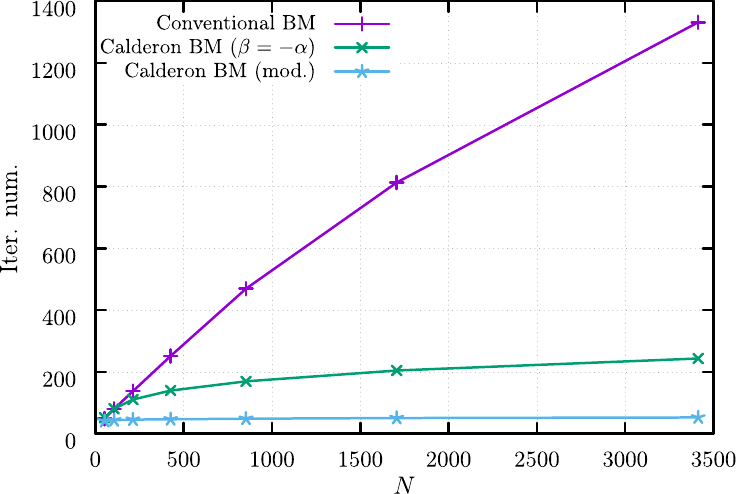}~~~
 \includegraphics[scale=0.35]{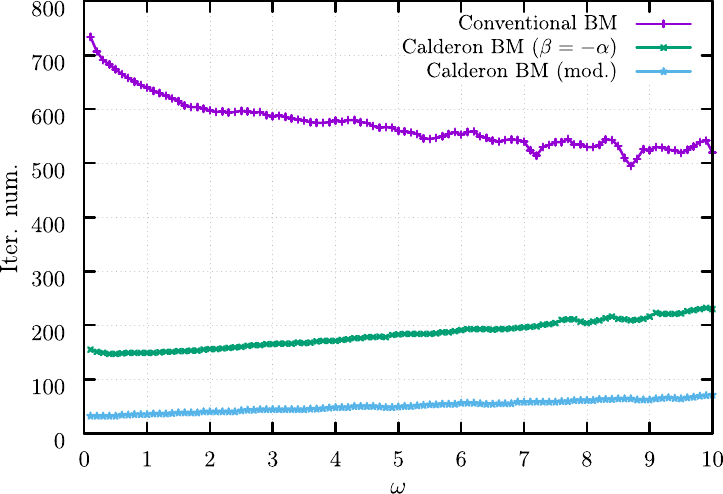}~~~
 \includegraphics[scale=0.35]{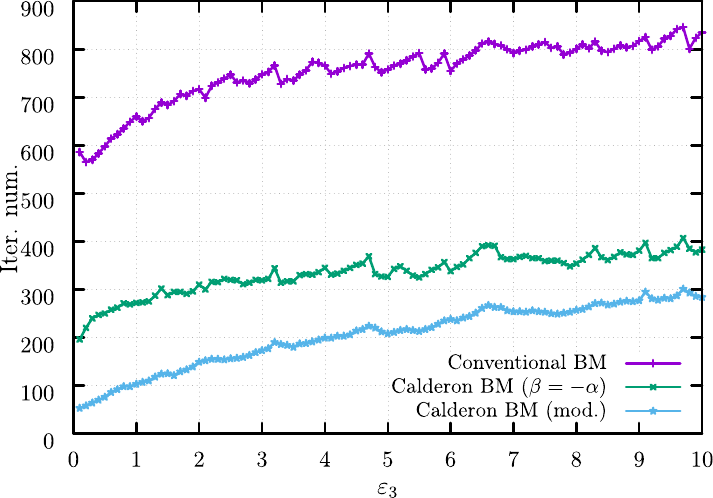}
 \caption{The number of GMRES iterations for the conventional and proposed integral formulations. The dependence on the total number of boundary elements $N$ (left), the angular frequency $\omega$ (centre), and the material constant $\varepsilon_3$ (right).}\label{fig:3mata}
\end{figure}

We then verify the performance of the proposed solvers in handling transmission problems defined in a more intricate domain, in which the boundaries of three different materials intersect at a single point (Figure \ref{fig:3mata_circle}). 
We here repeat the same computations for this setting as the previous concentric circular scatter case; evaluating the number of GMRES iterations by changing the number of boundary elements $N$, the angular frequency $\omega$, and the material constants $\varepsilon_2$ and $\varepsilon_3$. The values for $\omega$ and $\varepsilon_{2,\, 3}$ and are exactly the same as those for the previous examples. Figure \ref{fig:3mata} summarises the results. In the computations, the circular and vertical boundaries are respectively discretised into $3N/4$ and $N/4$ elements, and $N$ is fixed as $N=1066$ for centre and right figures.
One can again observe that the modified Calderon-preconditioned equation applying the Burton--Miller equation to some of the interior domains shows the fastest convergence for all the cases, followed by the original Calderon-preconditioned and then the conventional Burton--Miller ones. One may also notice that, compared with the results for previous concentric circular cases, the performance gap between the modified and original Calderon preconditioning has notably widened. This might be attributed to the mathematical and/or physical conditions of the underlying problems. The proposed formulation \eqref{eq:3dom_calderon_bie} may be insensitive to such conditions, which makes it more attractive.

\begin{figure}[h]
\centering 
\includegraphics[scale=0.8]{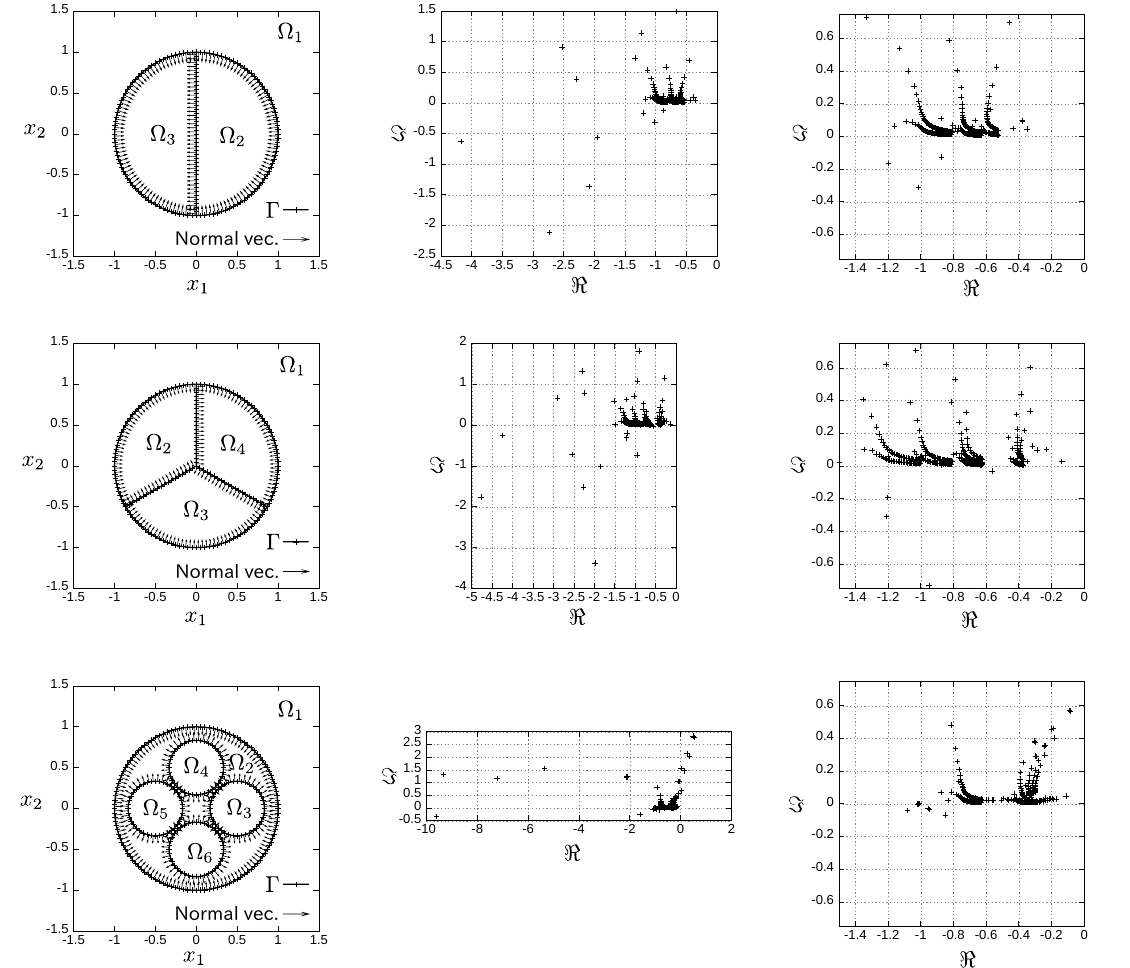}
\caption{The eigenvalue distribution of the squared coefficient matrices of the proposed boundary integral formulations. The left figures indicate the tested geometries as well as the normal on the boundaries. The centre figures show the whole eigenvalue distributions and the left ones are the magnified view of the centres.}
\label{fig:20_general_cases_eigenvalues}
\end{figure}

We finally check the eigenvalue distribution of the square of the proposed coefficient matrices for the previous test case (Figure \ref{fig:3mata_circle}) and some more general cases. The matrices are composed according to the recipe in Section \ref{sec:multi-material_cases}. The boundary elements of the tested geometries are shown in the left figures of Figure \ref{fig:20_general_cases_eigenvalues}. In the figures, the normal vectors are also plotted. The material constants for $\Omega_i$ indicated in the figures are given here as $\varepsilon_i=i$, and the angular frequency is set as $\omega=1.0$. The centre figures show all the eigenvalues of the squared coefficient matrices, and the right ones are their magnified versions around the clustering points. One finds that in all cases tested here, no eigenvalue has a significantly large magnitude. Also, there are only a few clustering points and they all are away from the origin. With these results, we conclude that the proposed integral equations are promising in efficiently solving the transmission problem defined in geometrically complicated domains.

\section{Conclusions}\label{sec:conclusions}
In this paper, we proposed novel boundary integral formulations based on the Burton--Miller formulation tailored for the transmission problem. We demonstrated that appropriately ordering the integral operators and multiplying proper constants for some rows and/or replacing some of the operators by the Burton--Miller ones can provide a well-conditioned integral formalism. Specifically, the square of the whole integral operator can be made equal to the one with a few eigenvalue \textcolor{black}{clustrering points}, providing a fast convergence for GMRES applied to the algebraic equations obtained by discretising the integral equation \textcolor{black}{via collocation}. Our formulation is easy to implement and can systematically be applied to various problems such as scattering by an object consisting of several materials with ``material junction points''.

It may be an interesting future direction to extend our formulation for transmission problems involving some vectorial waves such as electromagnetic~\cite{niino2012calderon} and/or elastic~\cite{isakari2012calderona} waves in three-dimensional space. It would also be interesting to apply other Krylov iterative solvers than GMRES such as some variants of BiCG and IDR$(s)$~\cite{sonneveld2009idr,sleijpen2010exploiting} to solve the algebraic equations stemming from the proposed integral equations. 

% \section*{Author contributions}
% Y. M. conceptualised the study, specifically developed the idea of rearranging the integral operators and found the best choice of the parameter $\beta$ for single-material case. A. Y. implemented test codes, corrected data, and drafted the manuscript. H. I. also conceptualised the study, specifically modifying Y. M.'s original formulation with the Calderon formulae and developed the multi-material formulations, and contributed to drafting and editing the manuscript while supervising the study.

\section*{Acknowledgments}
\textcolor{black}{The authors benefited from discussions with Mr. K. Tomoyasu from Keio University on the multi-material formulations.} This work was supported by JSPS KAKENHI Grant Numbers 23H03413, 23K19972 \textcolor{black}{and 24K20783}. \textcolor{black}{This work used the computational resources of the supercomputer Camphor3 provided by ACCMS of Kyoto University through the Joint Usage/Research Center for Interdisciplinary Large-scale Information Infrastructures and High Performance Computing Infrastructure in Japan (Project IDs: jh240031).}

\bibliographystyle{unsrt}
\bibliography{ref}

\end{document}